\documentclass[12pt]{amsart}

\usepackage{amssymb}
\usepackage{mathrsfs}

\def\F{\mathbb{F}}

\def\Z{\mathbb{Z}}

\newcommand{\pth}[1]{{#1}^{[p]}}

\def\su{{\subseteq}}
\def\la{{\langle}}
\def\ra{{\rangle}}

\def\w{{\omega}}

\def\dim{{\rm dim~}}
\def\char{{\rm char}}

\def\ad{\mbox{ad~}}

\def\Proof{\noindent{\sl Proof.}\ }
\def\qed{{\hfill $\Box$ \medbreak}}

\newtheorem{defi}{Definition}[section]
\newtheorem{thm}[defi]{Theorem}
\newtheorem{lem}[defi]{Lemma}
\newtheorem{cor}[defi]{Corollary}

\newtheorem{prob}[defi]{Problem}

\newtheorem{prop}[defi]{Proposition}

\begin{document}

\title[Polynomial identity enveloping algebras]{Enveloping algebras of restricted Lie superalgebras satisfying non-matrix polynomial identities}
\author{\textsc{Hamid Usefi}}
\thanks{ Research was supported by an  NSERC PDF}

\address{Department of Mathematics, University of British Columbia, 1984 Mathematics Road,
Vancouver, BC, Canada, V6T 1Z2}

\email{usefi@math.ubc.ca}

\begin{abstract}
Let $L$ be a restricted Lie superalgebra with  its enveloping algebra $u(L)$ over a field $\F$ of characteristic $p>2$.
A polynomial identity is called non-matrix if it is not satisfied by the algebra of $2\times 2$ matrices over $\F$.
We characterize $L$  when $u(L)$ satisfies a non-matrix polynomial identity.
In particular, we characterize $L$ when $u(L)$ is Lie solvable, Lie nilpotent, or Lie super-nilpotent.
\end{abstract}


  \maketitle

\section{Introduction}
A variety of associative algebras over a field $\F$ is called non-matrix if it does not contain $M_2(\F)$, the algebra of $2\times 2$ matrices over $\F$. A polynomial identity (PI) is called  non-matrix if $M_2(\F)$ does not satisfy this identity.
Latyshev in his attempt to solve the Specht problem proved that any non-matrix variety
generated by a finitely generated algebra over a field of characteristic zero is finitely based \cite{La1}.
The complete solution of the Specht problem in the case of characteristic zero is given by Kemer \cite{Kemer2, Kemer3}.

Although several counterexamples are found for the Specht problem in the positive characteristic \cite{AK},  the development in this area has lead to some interesting results. Kemer has investigated the relation between PI-algebras and nil algebras.
Amitsur \cite{Am} had already proved that the Jacobson radical of a relatively-free algebra of countable rank is nil.
Restricting to non-matrix varieties, Kemer \cite{Kemer1} proved that  the Jacobson radical of a relatively-free algebra of a non-matrix variety over a field of positive characteristic is nil of bounded index.
Recently these varieties have been further studied in \cite{MPR}.

Enveloping algebras satisfying polynomial identities were first considered by Latyshev \cite{L63} by proving that
the universal enveloping algebra of a Lie algebra $L$  over a field of characteristic zero satisfies a PI if and only if $L$ is abelian. Latyshev's result was extended to positive characteristic by Bahturin \cite{B74}.
Passman \cite{P90} and Petrogradsky \cite{P91} considered the analogous problem for restricted Lie algebras and their envelops.

Let $L=L_0\oplus L_1$ be a restricted Lie superalgebra with the bracket $(\,,)$.
 We denote the enveloping algebra of $L$ by $u(L)$. All algebras in this paper are over a field $\F$ of characteristic $p>2$ unless otherwise stated. In case $p=3$ we add the axiom $((y,y),y)=0$, for every $y\in L_1$.
 This identity is necessary to embed $L$ in $u(L)$.
 Restricted Lie superalgebras whose enveloping algebras satisfy a polynomial identity have been characterized by Petrogradsky \cite{P92}. The purpose of this paper is to characterize restricted Lie superalgebras whose enveloping algebras satisfy a non-matrix PI. Our results unify the results of Riley, Shalev, and Wilson in \cite{RS93, RW99} where they characterize restricted Lie algebras whose enveloping algebras satisfy a non-matrix PI. Our first main result is as follows.
 
\begin{thm}\label{uL-non-matrix}
The following statements are equivalent:
\begin{enumerate}
\item $u(L)$ satisfies a non-matrix PI.
\item The commutator ideal of $u(L)$ is nil of bounded index.
\item $u(L)$ satisfies a PI, $(L_0, L_0)$ is $p$-nilpotent, and there exists an  $L_0$-module $M$ of codimension at most 1 in  $L_1$ such that   $(M, L_1)$ is $p$-nilpotent and  $(L_1, L_0)\su M$.
\end{enumerate}
\end{thm}
We prove  Theorem \ref{uL-non-matrix} in  Section \ref{sec-nonmatrix}. The classification of finite dimensional Clifford algebras is used in the proof. In the course of proving Theorem \ref{uL-non-matrix} it was of interest for us to know whether a variant of Cayley-Hamilton Theorem holds for matrices over the Grassmann algebra. Note that the Grassmann algebra satisfies the identity $[x,y,z]=0$. So in general we ask the following:
\begin{prob}
Let $G$ be a nilpotent (solvable) Lie algebra over a field of positive characteristic. Does a
variant of Cayley-Hamilton Theorem hold for $M_n(G)$?
\end{prob}

Szigeti \cite{Sz} proved that if $H$ is a nilpotent Lie ring of any characteristic and $T\in M_n(H)$, then
there exists a polynomial $f(t)\in H[t]$ such that the left substitution of $T$ in $f(t)$ is zero. The degree of $f$ is $n^m$ where $m$ is the nilpotency class of $H$.  However the leading coefficient of $f$ is factorial in $n$, so this result is essentially not useful in  positive characteristic.

Every $\Z_2$-graded associative algebra $A=A_0\oplus A_1$ over $\F$ can be given the structure of a
restricted Lie superalgebra via $(x,y)=xy-(-1)^{ij}yx$ for every $x\in A_i$ and $y\in A_j$. We denote the usual Lie bracket of $A$ by $[u,v]=uv-vu$. We emphasis that the terms Lie nilpotent or Lie solvable are used with respect to the usual Lie bracket $[\, ,]$ whereas Lie super-nilpotent refers to the super-bracket $(\, ,)$.

The variety of Lie solvable algebras includes  Lie nilpotent  and Lie super-nilpotent algebras.
We characterize $L$ when $u(L)$ is Lie solvable in the following theorem.

\begin{thm}\label{uL-solvable}
Let $L$ be a restricted Lie superalgebra.
Then  $u(L)$ is Lie solvable if and only if
$(L, L)$ is finite-dimensional, $(L_0, L_0)$ is $p$-nilpotent, and there exists a subspace $M\su L_1$ of codimension at most 1 such that $(M, L_1)$ is $p$-nilpotent and $(L_1, L_0)\su M$.
\end{thm}
 Proof of Theorem \ref{uL-solvable} is given in Section  \ref{sec-solvable}. Furthermore, restricted Lie superalgebras whose enveloping algebras are   Lie super-nilpotent or Lie nilpotent are characterized in Theorem
 \ref{uL-Lie super-nilpotent} and Theorem \ref{uL-Lie nilpotent}, respectively.
Stewart \cite{St} proved that if $H$ is a nilpotent ideal of a Lie algebra $L$ such that $L/H'$ is nilpotent then
$L$ is nilpotent. In contrast we prove that if $I$ a nilpotent two-sided ideal of an associative algebra $R$
such that $R/I^2$ is Lie nilpotent then $R$ is Lie nilpotent, see Proposition \ref{R/I nilpotent}.

\section{Preliminaries}

Let  $A=A_0\oplus A_1$ be a vector space decomposition of a non-associative algebra. We
say that this is a $\mathbb{Z}_2$-grading of $A$ if $A_iA_j\su A_{i+j}$, for every
$i,j\in\Z_2$ with the understanding that the addition $i+j$ is mod 2.
The components $A_0$ and $A_1$ are called even and odd parts of $A$, respectively.
Note that $A_0$ is a subalgebra of $A$. One can associate a Lie super-bracket to $A$ by defining
$(x,y)=xy-(-1)^{ij}yx$ for every $x\in A_i$ and $y\in A_j$.

If $A$ is associative, then
for any $x\in A_i$, $y\in A_j$ and $z\in A$ the following identities hold:
\begin{enumerate}
\item[(1)] $(x,y)=-(-1)^{ij}(y,x)$,
\item[(2)] $(x,(y,z))=((x,y),z)+(-1)^{ij}(y,(x,z)).$
\end{enumerate}
The above identities are the defining relations of Lie superalgebras.
Furthermore, $A$ can be viewed as a Lie algebra by the usual Lie bracket
$[u,v]=uv-vu$.  Let  $B$ and $C$ be subspaces of $A$.
We denote by $[B, C]$ the subspace spanned by all commutators $[b, c]$, where $b\in B$ and $c\in C$.
Then $[B,_n C]=[[B,_{n-1} C], C]$, for every $n\geq 2$.
If $B$ is a nilpotent subalgebra of $A$, the nilpotence index of $B$ is the least integer $k$ such
that $B^{k+1}=0$. The lower Lie central series of $A$ is defined  by setting
 $\gamma_1(A)=A$ and  $\gamma_n(A)=[\gamma_{n-1}(A), A]$, for every $n\geq 2$.
  The Lie derived series of $A$ is defined by  setting $\delta_0(A)=A$ and   $\delta_n(A)=[\delta_{n-1}(A), \delta_{n-1}(A)]$, for every $n\geq 1$.
  Recall that $A$ is called Lie nilpotent if $\gamma_n(A)=0$, for some $n$, and
$A$ is called Lie solvable if $\delta_m(A)=0$, for some $m$.
The nilpotence class of $A$ is the least integer $c$ such that $\gamma_{c+1}(A)=0$.
The derived length of $A$ is the least integer $d$ such that $\delta_d(A)=0$.
Long commutators are left tapped, that is $[x,y,x]=[[x,y], z]$, and $[x,_k y]$ denotes the commutator $[x, y, \ldots, y]$, where $y$ occurs $k$ times. Analogous definitions hold for the super-bracket of $A$.
We denote by $\gamma_n^s(A)$  the $n$-th term of the Lie super-central series of $A$.

If $L$ is a Lie superalgebra we denote the bracket of $L$ by $( , )$.
The adjoint representation of $L$ is given by $\ad x : L\to L$, $\ad x(y) = (y, x )$, for all $x,y\in L$.
The notion of restricted Lie superalgebras can be easily formulated as follows:
\begin{defi}\label{def:res}
A Lie superalgebra $L = L_0\oplus L_1$ is called  \em
{restricted}, if there is a $p$th power map
$L_0\to L_0$, denoted as $\pth{}$, satisfying
\begin{enumerate}
\item[(a)] $\pth{(\alpha x)} = \alpha^p \pth x$, for all  $x \in L_0$ and $\alpha \in \F$,
\item[(b)] $(y, \pth x) = (y,_p x)$, for all $x \in L_0$
and $y \in L$,
\item[(c)] $\pth{(x + y)} = \pth x + \pth y + \sum_{i= 1}^{p-1}
s_i(x,y)$, for all $x, y \in L_0$ where $is_i$ is the
coefficient of ${\lambda}^{i-1}$ in $(\ad(\lambda x +y))^{p-1}(x)$.
\end{enumerate}
\end{defi}
In short, a restricted Lie superalgebra is a Lie superalgebra
whose even subalgebra is a restricted Lie algebra and the odd part
is a restricted module by the adjoint action of the even
subalgebra. For example, every $\mathbb{Z}_2$-graded associative algebra inherits a restricted Lie
superalgebra structure.

Let $L$ be a restricted Lie superalgebra. We denote the enveloping algebra of $L$ by $u(L)$.
The augmentation ideal $\w(L)$ is the ideal of $u(L)$ generated by $L$.
The analogue of the PBW Theorem is as follows. We refer to \cite{B92} for basic background.

\begin{thm}
Let $L=L_0\oplus L_1$ be a restricted  Lie superalgebra and let ${\mathcal B}$ be a totally ordered
basis for $L$ consisting of $\Z_2$-homogeneous elements. Then $u(L)$ has a basis
consisting of PBW monomials, that is, monomials of the form $x_{1}^{a_1}\ldots x_{s}^{a_s}$ where
$x_1<\cdots < x_s$ in ${\mathcal B}$, $0\leq a_i<p$ whenever $x_i\in L_0$, and $0\leq a_i\leq 1$ whenever $x_i\in L_1$.
\end{thm}

The fact that $u(L)$ is $\Z_2$-graded can be explained as follows.
For every Lie superalgebra $L=L_0\oplus L_1$ we can associate a linear map
$\sigma: L\to L$ such that $\sigma^2=1$. Indeed, $\sigma(x)=x$, for all $x\in L_0$ and $\sigma(y)=-y$, for all $y\in L_1$.
The converse of this is also true, that is every linear map of $L$ of order 2 induces a $\Z_2$-grading on $L$.
Now suppose that $L=L_0\oplus L_1$ is a Lie superalgebra and let $\sigma$ be the corresponding linear map.
Then $\sigma$ extends to an automorphism of $A=u(L)$ of order 2. Now, we can write
$A=A_0\oplus A_1$, where
$$
A_0=\la u\in A\mid \sigma(u)=u\ra_{\F}, \quad \text{and} \quad A_1=\la u\in A\mid \sigma(u)=-u\ra_{\F}.
$$
Since $\sigma$ is an automorphism of $A$, the parity of a PBW monomial $x_{1}^{a_1}\ldots x_{s}^{a_s}$
is equal to the parity of the number of odd $x_i$ with exponent 1.

Since $L$ embeds into $u(L)$ we identify $\pth x$ with $x^p$, for every $x\in L_0$.
Furthermore, if $x\in L_0$ and $y\in L$ then the bracket $(x,y)$ in $L$ is the same as the bracket $[x, y]$ in $u(L)$.
By an ideal of $L$ we always mean a restricted ideal, that is $I$ is an ideal of $L$ if
 $(I, L)\su I$ and $I_0$ is closed under the $p$-map.
Let $H$ be a subalgebra of $L$. We denote by $H'$ the commutator subalgebra of $H$, that is $H'=(H, H)$.
For a subset $X\su  L$, we denote by $\la X\ra_p$ the restricted ideal of $L$ generated by $X$. Also,
 $\la X\ra_{\F}$ denotes the subspace spanned by $X$.
 An element $x\in L_0$ is called \emph{$p$-nilpotent} if there exists some non-negative integer $t$ such that $x^{p^t}=0$. Also, recall that $X$ is said to be $p$-nil if every element $x\in X$ is $p$-nilpotent and  $X$ is $p$-nilpotent if there exists a positive integer $k$ such that $x^{p^k}=0$, for every $x\in X$.

We shall call $L$  a nilpotent $L_0$-module if $(L, _n L_0)=0$, for some $n$.
Note that Engel's Theorem holds for Lie superalgebras, see \cite{sch}, for example.
\begin{thm}[Engel's Theorem]
Let $L$ be a finite-dimensional Lie superalgebra such that
$\text{ad } x$ is nilpotent, for every homogeneous element $x\in L$.
Then $L$ is nilpotent.
\end{thm}

The proof of the following lemma follows from Engel's Theorem and the fact that $(\text{ad } x)^2=\frac{1}{2} \text{ad } (x, x)$, for every $x\in L_1$.

\begin{lem}\label{L-res-nilpotent}
Let $L=L_0\oplus L_1$ be a finite-dimensional restricted Lie superalgebra.
If $L_0$ is  $p$-nilpotent then $L$ is nilpotent.
\end{lem}

\begin{lem}\label{w(L)-nilpotent}
Let $L$ be a restricted Lie superalgebra. Then $\w(L)$ is associative nilpotent if and only if
$L$ is finite-dimensional  and $L_0$ is $p$-nilpotent.
\end{lem}
\Proof The if part  follows from the PBW Theorem.
We prove the converse by induction on $\dim L$. Since $L$ is nilpotent, by Lemma \ref{L-res-nilpotent},
there exists a non-zero element $z$ in the center $Z(L)$ of $L$. Since $Z(L)$ is homogeneous we may assume that
either $z\in L_0$ or $z\in L_1$. If $z\in L_1$ then $z^2=(z,z)=0$. If $z\in L_0$ then
we can replace $z$ with its $p$-powers so that $z^p=0$. So in either case $z^p=0$ in $u(L)$.
Now consider $H=L/\la z\ra_p$. Then by induction hypothesis $\w(H)$ is nilpotent.
This means that $\w^m(L)\su \la z\ra_p u(L)$, for some $m$. It then follows that
$\w^{mp}(L)\su \la z^p\ra_p u(L)=0$, as required.
\qed

We shall need the following result of Petrogradsky.

\begin{thm}[\cite{P92}]\label{Petrogradsky}
Let $L=L_0\oplus L_1$ be a restricted Lie superalgebra.
Then $u(L)$ satisfies a  PI if and only if there exist homogeneous restricted  ideals $M\su N\su L$ such that
\begin{enumerate}
\item $\dim L/N$ and $\dim M$ are both finite.
\item $N'\su M$, $M'=0$.
\item The restricted Lie subalgebra  $M_0$ is $p$-nilpotent.
 \end{enumerate}
\end{thm}

The Grassmann (or exterior) algebra $G$ on a vector space $V$ is defined by the relation:
 $v^2=0$, for all $v\in V$. This relation implies that $uv = - vu$, for all $u,v\in V$.
Note that $G$ satisfies $[x,y,z]=0$. Since  the identity $(x+y)^p=x^p+y^p$ holds in any Lie nilpotent algebra of class at most $p$, it follows that $G$ is nil of index $p$.

Let $V$ be an $\F$-vector space and $Q:V\to \F$ a quadratic form, that is $Q$ satisfies the following conditions:
\begin{itemize}
\item[(1)] $Q(\alpha x)=\alpha^2 Q(x), x\in V, \alpha\in \F$.
\item[(2)] $B(x,y)=Q(x+y)-Q(x)-Q(y)$ is bilinear.
\end{itemize}
 Let $T(V)=\F\oplus V\oplus V\otimes V\oplus\cdots$ be the tensor algebra on $V$. The associated Clifford algebra
 to $V$ and $Q$ is defined by $C(V, Q)=T(V)/\la x\otimes x-Q(x)1\ra$. It is known that if  $Q$ is non-degenerate and
 $\dim V=n$ is even then $C(V,Q)$ is a central simple algebra of dimension $2^n$, see \cite{Ja} for example.

\section{General non-matrix PI}\label{sec-nonmatrix}

Let $R$ be an algebra satisfying a non-matrix polynomial identity $f$. It follows from the structure theory of PI-algebras that
$[R, R]R$ is nil. Indeed, by Posner's theorem, if a prime ring  satisfies $f$, it must be commutative.  Since every semiprime ring is a subdirect sum of prime rings, if a semiprime ring satisfies $f$, it must be commutative.  But $R/N$ is semiprime,
where  $N$ is the nil radical of $R$ (the sum of all nil two-sided ideals). So, $R/N$ must be commutative.
Thus, $[R, R]R$ is nil. However, a stronger result holds using a theorem of Kemer.

\begin{thm}[\cite{Kemer1}]\label{kemer}
The Jacobson radical of a relatively-free algebra of a non-matrix variety over a field of positive characteristic is nil of bounded index.
\end{thm}

\begin{cor}\label{ass-PI}
If $R$ satisfies a non-matrix PI over a field of positive characteristic, then  $[R,R]R$ is nil of bounded index.
\end{cor}
\Proof By Theorem \ref{kemer}, the radical $J$ of the relatively-free algebra $U$ of the variety defined by $R$ is nil of bounded index. We note that $U/J$ is commutative since it is semisimple satisfying a non-matrix PI. Thus,
$[U,U]U\su J$, implying that $[U, U]U$ is nil of bounded index.  Thus $[R,R]R$ is nil of bounded index.
\qed

\begin{prop}\label{M}
If $u(L)$ satisfies a non-matrix PI then there exists an $L_0$-module  $M$  of codimension at most 1 in $L_1$
such that $(M, L_1)$ is $p$-nilpotent and $(L_1, L_0)\su M$.
\end{prop}
\Proof Let $M$ be the set of all $y\in L_1$ such that $(y,y)$ is $p$-nilpotent.
By Theorem \ref{Petrogradsky}, there exists a homogeneous ideal $N$ of $L$ of finite codimension
such that  $(N_1, N_1)$ is $p$-nilpotent. In particular, $N_1\su M$.
Since $[R,R]R$ is nil of bounded index, by Corollary \ref{ass-PI}, it follows that $(M, L_1)$ is $p$-nilpotent.
Note that $(L_1, L_0)\su M$. It is enough to show that the codimension of $M$ in $L_1$ is at most 1.
Suppose to the contrary. Let $K=\la (L_0, L_0) + M\ra_p$.
Without loss of generality we can replace $L$ with $L/M$.
So,  $(a,a)$ is not $p$-nilpotent, for every $a\in L_1$.
Let $a$ and $b$ be linearly independent elements of $L_1$
and set $x=(a,a)$, $y=(b,b)$, and $z=(a,b)$.
We replace $L$ with its subalgebra $H=L_0\oplus \la a,b\ra_{\F}$.
Note that $(a,b)$ is not $p$-nilpotent.
Let $X$ be the subset of $R$ consisting of all $x^iy^jz^k$, where $i,j,k$ are non-negative integers. Note that $X$ is a multiplicative set,  $1\in X$ and $0\notin X$.
Now we consider the localization  $R_{_X}$ of $R$ with respect to $X$.
Consider the subalgebra $S$ of $R_{_X}$ generated by all $r/s$, $r,s\in X$.
Let $\mathfrak{m}$ be a maximal ideal of $S$ and set $F=S/\mathfrak{m}$.
Let $V=H_1\otimes_{\F} F$. We denote by $\bar r$ the cost representative of an element $r\in S$.
Consider the $F$-bilinear map $\phi: V\times V\to F$ induced by $\phi(u,v)=\overline{(u,v)}$, for every $u,v\in H_1$.
Then the Clifford algebra $C$ over $F$ defined by $\phi$ is the $F$-subalgebra of $R_F=R_{_X}\otimes_{\F} F$ generated by $a,b$.
By the classification of Clifford algebras, $C$ is a central simple algebra.
We also claim that  $R_F$ satisfies a non-matrix PI. Note that $(H_0, H_1)=0$. We observe that  the commutator ideal $[R_F, R_F]R_F$ of $R_F$ is generated by $[a,b]\otimes 1,  [a,b]a\otimes 1$, and   $[a,b]b\otimes 1$.
Since each of these generators are nil, by Corollary \ref{ass-PI}, it follows that $[R_F, R_F]R_F$ is nilpotent.
 Thus, $R_F$ is nilpotent-by-commutative which is a non-matrix PI.
Thus, $C$ satisfies a non-matrix PI and being central simple,  $C$ must be commutative. So, $ab = ba$.  But then $2ab=(a,b) \in L_0$ which contradicts the PBW Theorem because $a$ and $b$ are linearly independent. This contradiction implies that $L_1$ is 1-dimensional. Hence, $M$ has codimension at most 1 in $L_1$.\qed

\begin{lem}\label{ideal-nil}
Let $L$ be a restricted Lie superalgebra and $N$  a homogeneous abelian ideal of $L$ of finite codimension.
Let $I$ be an ideal of $u(N)$ that is stable under the adjoint action of $L$ on $N$.
If $I$ is nil of bounded index then so is  $Iu(L)$.
\end{lem}
\Proof Let $R=u(L)$.
Let $x_1,\ldots, x_n\in L_0$ and $y_1, \ldots, y_m\in L_1$ such that
$$
L=N+\la x_1,\ldots, x_n, y_1, \ldots, y_m\ra_{\F}.
$$
By the PBW Theorem $R$ has a basis consisting of the monomials of the form
\begin{align*}
&wx_1^{\alpha_1}\ldots x_n^{\alpha_n}y_1^{\beta_1}\ldots y_m^{\beta_m}\\
&0\leq \alpha_i<p, \quad \beta_j\in \{0,1\},
\end{align*}
where the $w$'s are PBW monomials lying in $u(N)$.
Let $D=u(N)$. Note that $R$ is a left $D$-module of finite rank $r=p^n2^m$.
Thus, under the left regular representation of $R$, we can embed $R$ into $M_r(D)$.
Since $I$ is stable under the adjoint action of $L$ on $N$, we have
$RI=IR=RIR$. We observe that $IR$ embeds in $M_r(I)$.
Therefore, it suffices to show that $M_r(I)$ is nil of bounded index.
Note that $u(N)$ is the tensor product of a commutative algebra with the Grassmann algebra.
Thus, $u(N)$ satisfies the identity $[x, y, z]=0$. Hence, $I$ satisfies $[x, y, z]=0$, too.
Recall  Levitzki's Theorem and Shirshov's Height Theorem
stating that every $t$-generated PI algebra which is nil of bounded index $s$ is (associative) nilpotent of a bound given as a function of $s$, $t$, and $d$, where $d$ is the degree of the polynomial identity, see \cite{Lev} and \cite{Shirshov}.
So, if $S$ is any $r^2$-generated subalgebra of $I$ then there exists a constant $k$ such that
 $S^k=0$. Now, let $T\in M_r(I)$ and denote by $S$ the subalgebra of $I$ generated by all entries of $T$.
So,  $T^i\in M_n(S^i)$, for every $i$. Since $S^k=0$, we get $T^k=0$. Since $k$ is independent of $T$,
it follows that $M_r(I)$ is nil of bounded index, as required.
\qed

\noindent\emph{Proof of Theorem \ref{uL-non-matrix}.}
  $(1)\Leftrightarrow (2)$: follows from  Corollary \ref{ass-PI}.\\
$(1)\Rightarrow (3)$:  follows from  Proposition \ref{M} and Corollary \ref{ass-PI}.\\
$(3)\Rightarrow (2)$: we shall use the fact that the class of   nil algebras of bounded index is closed under extensions.
Let $L_1=M+\F z$, where $(M, L_1)$ is $p$-nilpotent and $(L_1, L_0)\su M$.
Since $u(L)$ satisfies a PI, by Theorem \ref{Petrogradsky}, there exists a homogenous restricted ideal $N$ of $L$ such that $N'$ is finite-dimensional and $(N')_0$ is  $p$-nilpotent.
Note that $(N, N)$ is nilpotent, by Lemma \ref{L-res-nilpotent}.
Thus, $\w(N')$ is associative nilpotent, by Lemma \ref{w(L)-nilpotent}. Hence, $N'u(L)$ is associative  nilpotent.
Consider the natural map
$u(L)\to u(L/N')$.
It suffuses to consider  $L/N'$. In other words, we may assume
that $N$ is abelian. Let $B$ be the restricted ideal of $L$ generated by $(N, L)$.
We observe that  $I=Bu(N)$ is nil of bounded index.
Furthermore, by Lemma \ref{ideal-nil}, $Iu(L)=Bu(L)$ is also a nil ideal of $R$ of bounded index.
But $Iu(L)$ is the kernel of the homomorphism $R\to u(L/B)$. Thus, we can replace
$L$ with $L/B$. In other words, we can assume that $N$ is central in $L$.
It follows that $L'$ is finite-dimensional. Note that $N_1u(N)$ is  nil of index $p$. Since $N_1$ is a central ideal of $L$, $N_1R$ is nil of bounded index, by Lemma \ref{ideal-nil}. Thus, we may assume that $N_1=0$.
It follows that $M$ is finite-dimensional. Thus, by Lemma \ref{w(L)-nilpotent}, $((M, L_1)+M)u(L)$ is associative  nilpotent. We can now assume $M=0$. Hence, $L_1=\F z$ and it follows that $(z, L_0)=0$.
Note that $(zR)^2=0$. Therefore, we can assume that $L$ is a restricted Lie algebra for which $L'$ finite-dimensional  and $p$-nilpotent. Thus, the associative ideal of $R$ generated by $L'$ is associative nilpotent, by Lemma \ref{w(L)-nilpotent}. So, we may assume that $L'=0$. But then $R$ is a commutative algebra and the assertion holds.
\qed

\section{Lie Solvable}\label{sec-solvable}

\begin{prop}\label{derived-sub-f.d.}
 If $u(L)$ is Lie solvable then $\dim (L, L)$ is finite.
\end{prop}

The following result is proved by Zalesskii and Smirnov \cite{ZS} and independently by Sharma and
Srivastava \cite{SS}.

\begin{thm}\label{SS}
Let $R$ be a Lie solvable ring  of derived length $t\geq 2$. The two-sided ideal $\mathscr{J}$ of $R$
generated by  $[[R, R], [R, R], R]$ is associative nilpotent of index bounded by a function of $t$.
\end{thm}

Throughout this section, we denote $u(L)$ by $R$ and  $\mathscr{J}$ is used to denote the associative ideal of $u(L)$ generated by
 $[[R, R], [R, R], R]$.

The proof of Proposition \ref{derived-sub-f.d.} breaks down to several parts.
First we prove that $(L_0, L_0)$ is finite-dimensional.
In this case, it is enough to assume that $L$ is a  restricted Lie algebra.
We remark that  this assertion is proved in \cite{RS93}, however we can offer a new shorter proof.
Let us first recall the following result of Neumann originating from Group theory, see \cite{Neumann, B}.

\begin{thm}\label{Neumann}
Let $\phi: U\times V \to W$ be a bilinear map, where $U, V, W$ are
vector spaces over $\F$. Suppose that there exist constants $m,\ell$ such that
 for each $u\in U$ the codimension of its annihilator in $V$ is bounded by $m$, and for each $v\in V$ the codimension of its annihilator in $U$ is bounded by $\ell$. Then $\dim \la \phi(U, V)\ra_{\F} \leq m\ell$.
\end{thm}

\begin{lem}\label{derived-even-f.d.} Let $K$ be a restricted Lie algebra.
If $u(K)$ is Lie solvable then $K'$ is finite-dimensional.
\end{lem}
\Proof Clearly $K$ is solvable. \\

\noindent\emph{Step 1.} We can assume $\delta_2(K)=0$.
We use induction on the derived length $s$ of $K$. If $K'=0$, there is nothing to prove. So, assume that $s\geq 2$.
  Let $H=K/\la \delta_{s-1}(K)\ra_p$. Since $u(H)$ is Lie solvable, we can assume by induction hypothesis that $[H, H]$ is finite-dimensional. Since $[K, K]$ is $p$-nilpotent, by Corollary \ref{ass-PI}, it is enough to prove that $\delta_{s-1}(K)$ is finite-dimensional. Without loss of generality we can replace $K$ with $\la \delta_{s-2}(K)\ra_p$.
So, $\delta_2(K)=0$.\\

\noindent\emph{Step 2.} We can assume $\gamma_3(K)=0$.
We prove that $\gamma_3(K)$ is finite-dimensional. Then we replace $K$ with $K/\la \gamma_3(K)\ra_p$.
We apply Theorem \ref{Neumann} to the function $\varphi:K\times K'\to \gamma_3(K)$ given by $\varphi(x,y)=[x,y]$, for every $x\in K$ and $y\in K'$.
It is enough to prove that $\dim [x, K']$  and $\dim [y, K]$ are bounded for every $x\in K$ and $y\in K'$.
Fix $x\in K$. For every $y\in K'$, we write $y'=[y, x]$.
Since $K'$ is abelian, we have
\begin{align*}
[y_1x, y_2x]+[y_1y_2x, x]=2y_1'y_2x\in [R, R],
\end{align*}
for every $y_1,y_2\in K'$.
Since $[yx, x]=y'x\in [R, R]$ and $\char(\F)\neq 2$, we have
\begin{align*}
[y_1'y_2x, y_1'x]=y_1'^2y_2'x\in \delta_2(R),
\end{align*}
Thus,
\begin{align*}
[y_1'^2y_2'x, y_2]=-y_1'^2y_2'^2\in  \mathscr{J},
\end{align*}
which implies that $\dim [x, K']$ is bounded since $\mathscr{J}$ is associative nilpotent.
Now let $y\in K'$ and set  $\dot x=[x, y]$, for every $x\in K$.
 We also define $x_1\circ x_2= x_1\dot x_2+ \dot x_1 x_2$, for every $x_1,x_2\in K$.
Then $x_1\circ x_2=[x_1x_2, y]\in [R, R]$ and $[x_1x_2y, y]=(x_1\circ x_2)y\in [R, R]$.
Since $K$ is abelian, we have
$$
[x_1\circ x_2, (x_1\circ x_2)y]=2(x_1\circ x_2)\dot x_1\dot x_2\in \delta_2(R).
$$
Hence,
$$
 [(x_1\circ x_2)\dot x_1\dot x_2, y]=2\dot x_1^2 \dot x_2^2\in \mathscr{J},
 $$
 which implies that $\dim [y, K]$ is bounded since $\mathscr{J}$ is associative nilpotent.\\

\noindent\emph{Step 3.} $K'$ is finite-dimensional.
We apply Theorem \ref{Neumann} to the function $\varphi:K\times K\to K'$ given by $\varphi(x,y)=[x,y]$, for every $x, y\in K$.
Thus, it suffices to show that $\dim [x, K]$ is  bounded for every $x\in K$. This is similar to Step 1 taking into account that
$\gamma_3(K)=0$. Thus, $K'$ is finite-dimensional, as required.
\qed

By Theorem \ref{Petrogradsky} there exists a homogenous restricted ideal $N$ of $L$ such that $N'$ is finite-dimensional.
In order to prove that $(L, L)$ is finite-dimensional it suffices to replace $L$ with $L/N'$.
In the next two lemmas, we assume that $L$ has an abelian ideal $N$ of finite codimension.

\begin{lem}\label{(L_1, L_1)}
If $u(L)$ is Lie solvable then $(L_1, L_1)$ is finite-dimensional.
\end{lem}
 \Proof We apply Theorem \ref{Neumann} to the  function $\varphi: L_1\times L_1\to (L_1, L_1)$ given by $\varphi(y, z)=(y, z)$, for $y, z\in L_1$.  Since $N$ is abelian and $\dim L/N$ is finite, it is enough to prove  that $\dim (y, N_1)$ is finite, for every $y\in L_1$.   Fix $y\in L_1$ and let $z'=(z,y)$, for every $z\in L$.
  Note that $[z_1z_2, z_3]=z_1(z_2,z_3)-(z_1, z_3)z_2$, for every  $z_1,z_2,z_3\in L_1$.
   Since $N$ is abelian, we have
   $[y_1, y_2]=2y_1y_2\in [R, R]$ and $[y_1, y]=y_1'-2yy_1\in [R, R]$, for every $y_1, y_2\in N_1$.
   Since $\char(\F)\neq 2$, we get
   $$
   [y_1'-2yy_1, y_1y_2]=2[y_1y_2, y]y_1=2y_1'y_1y_2\in \delta_2(R).
   $$
   Thus,
   $[y_1'y_1y_2, yy_2]=y_1'y_1y_2'y_2\in \mathscr{J}$.
We deduce that $\dim (y, N_1)$ is finite as $\mathscr{J}$ is nilpotent  by Theorem \ref{SS}.
  \qed

In order to finish the proof of Proposition \ref{derived-sub-f.d.}, it remains to prove the following:

\begin{lem}\label{even-odd}
If $u(L)$ is Lie solvable then $(L_1, L_0)$ is finite-dimensional.
 \end{lem}
 \Proof Consider the function $\varphi: L_0\times L_1\to (L_0, L_1)$ given by $\varphi(x,y)=(x,y)$, for
 every $x\in L_0$ and $y\in L_1$. By Theorem \ref{Neumann} and using the fact that  $N'=0$ and $\dim L/N$ is finite, it is enough to prove that $\dim (N_0, y)$ and   $\dim (N_1, x)$ are finite, for every $x\in L_0$ and $y\in L_1$. Since $N$ is abelian, we have
 \begin{align*}
 [[x_1, x_3y], [x_2, y], y]&=2[x_1, y][x_2,y][x_3,y]\in \mathscr{J},\\
 [[x, y_1], [x, xy_2], y_3]&=2[x, y_1][x, y_2][x,y_3]\in \mathscr{J},
\end{align*}
 for every $x_1, x_2, x_3\in N_0$ and $y_1, y_2, y_3\in N_1$. Since, by Theorem \ref{SS}, $\mathscr{J}$ is nilpotent, the assertion follows.
 \qed

We next show that if $R$ is Lie super-nilpotent then $R$ is in fact Lie solvable.

\begin{lem}\label{super-nilp-solvable}
Let $R=R_0\oplus R_1$ be any superalgebra. If $\gamma_{c+1}^s(R)=0$ then $\delta_c(R)\subseteq R_0$. Consequently, $R$ is Lie solvable.
\end{lem}
\Proof Observe that $\delta_1(R)=[R,R]\subseteq R_0+(R_1,R_0)$. Hence,
$$
\delta_2(R)\subseteq R_0+(R_1,R_0,R_0).
$$
Continuing this way yields $\delta_m(R)\subseteq R_0+(R_1,_{m} R_0)$, for every $m\geq 1$. So, $\delta_c(R)\subseteq R_0$.
\qed

\begin{thm}\label{uL-Lie super-nilpotent}
Let $L=L_0\oplus L_1$ be a restricted Lie superalgebra.
Then  $u(L)$ is Lie super-nilpotent if and only if $L$ is nilpotent,  $L'$ is finite-dimensional, $(L_0, L_0)$ is $p$-nilpotent, and  there exists an $L_0$-module $M$ of codimension at most 1 in  $L_1$ such that  $(M,L_1)$ is $p$-nilpotent and $(L_1, L_0)\su M$.
\end{thm}
\Proof Suppose that $R=u(L)$ is Lie super-nilpotent. Then, by Lemma \ref{super-nilp-solvable}, $R$ is  Lie solvable.
So, the if part follows from Corollary \ref{ass-PI} and Propositions \ref{M} and  \ref{derived-sub-f.d.}.
We prove the converse.  Let $z\in L_1$ such that $L_1=M+\F z$ and set $H=L_0+M$.\\

\noindent\emph{Step 1.} We may assume $(H, H)=0$.
 Since $H$ is nilpotent, there exists a non-zero element $x\in (H, H)\cap Z(H)$, where
 $Z(H)$ is the center of $H$.
Note that $(x, z)\in (H, H)\cap Z(H)$. Since $L$ is nilpotent, we can replace $x$ with $(x,_k z)$, for some $k$, so that $x\in Z(L)$. Since $(L,L)$ and $Z(L)$ are graded,
we can assume either $x=x_0$ or $x=x_1$. If $x\in L_0$ then either
$x\in (L_0, L_0)$ or $x\in (M, L_1)$. In either case, $x$ is $p$-nilpotent.
Let $s$ be the least integer such that $x^{p^{s+1}}= 0$. Then we replace $x$ with $x^{p^{s}}$.
Thus, $x^p=0$.
On the other hand, if $x\in L_1$ then $(x,x)=0$. Thus, in either case $x^p=0$.
Now $I=\la x\ra_{\F}$ is a restricted ideal of $L$ and we consider $L/I$.
Suppose that $u(L/I)=R/IR$ is Lie super-nilpotent. This means that there exists an integer $m$ such that
$
(R,_m R)\subseteq IR.
$
Since $(I, R)=0$, we get
$$
(IR,_{(p-1)m} R) \su  I^pR=0.
$$
Hence, $(R,_{pm} R)=0$. So it is enough to prove that $u(L/I)$ is Lie super-nilpotent.
Hence, by induction on  $\dim (H,H)$, we can assume $(H,H)=0$.\\

\noindent\emph{Step 2.} We may assume $(H, z)=0$.
Let $y$ be a homogeneous element in $H$ such that $x=(y, z)\neq 0$.
Since $x\in H$, we can replace $x$ with $(x,_i z)$, for some $i$,  so that $(x, z)=0$.
So, $x\in Z(L)\cap (L, L)$. Now we can use a similar argument as in Step 1, to show that we can replace $L$
with $L/\la x\ra_{\F}$. So, by induction on  $\dim (H, z)$, we can assume $(H, z)=0$.

By Steps 1 and 2, we can assume $(H, L)=0$. It is now easy to verify that $\gamma_3^s(R)=0$, as required.
\qed

\noindent\emph{Proof of Theorem \ref{uL-solvable}.}
The if part follows from Corollary \ref{ass-PI} and Propositions \ref{M} and  \ref{derived-sub-f.d.}.
We prove the converse.
Let $K$ be the ideal of $L$ generated by $(L, L_0)$. Since $K_0$ is $p$-nilpotent and
$(L, L)$ finite-dimensional, it follows by Lemma \ref{w(L)-nilpotent} that $\w(K)$ is associative
nilpotent. Thus, $Ku(L)$ is associative nilpotent. So, it is enough to prove that
$u(L/K)=u(L)/Ku(L)$ is Lie solvable. Hence, we can assume $L_0$ is central in $L$.
In a similar manner, we can assume $(M, L)=0$. Thus, $L$ is nilpotent and it follows from Theorem \ref{uL-Lie super-nilpotent} that $u(L)$ is Lie super-nilpotent.
Thus, by Lemma \ref{super-nilp-solvable}, $u(L)$ is Lie solvable, as required.
\qed

\section{Lie nilpotence}\label{sec-Lie nilpotence}

In this section we characterize $L$ when $u(L)$ is Lie nilpotent.

\begin{lem}\label{nilp-[L,L]-p-nilp}
If  $u(L)$ is Lie nilpotent then either $(L_1, L_1)$ is $p$-nilpotent or $\dim L_1\leq 1$ and $(L_1, L_0)=0$.
\end{lem}
\Proof Suppose that $(L_1, L_1)$ is not $p$-nilpotent. By Proposition  \ref{M}, there exists a subspace $M$ of $L_1$ of codimension at most 1 such that $(M, L_1)$ is $p$-nilpotent. Let  $L_1=M+\F z$. Then $(z,z)$ is not $p$-nilpotent. Let $y\in M$.
So,  $(z,z)$ and $(y,z)$ are linearly independent.
Now for every $n\geq 1$, we have
\begin{align}\label{yix}
[y,_{n} z]=(-2)^n z^ny+ (y,_{n} z)+ \sum_{i+j=n} \alpha_{i,j}z^j (y, _{i} z) ,
\end{align}
where in the sum above $i,j$ are positive integers and $\alpha_{i,j}\in \F$.  Let $m$ such that $[y,_{2m} z]=0$.  Note that
 $(z^2)^my$ is a PBW monomial of degree $m+1$ and the rest of the PBW monomials in Equation (\ref{yix}) have degree strictly less than $m+1$. Thus, $(z^2)^my=0$. Since $(z,z)$ is not nilpotent, we deduce that $y=0$. Thus, $M=0$ and $L_1$ is 1-dimensional.
Let $x\in L_0$. Then  $[z,x]=\alpha z$, for some $\alpha \in \F$. But $[z,{_n} x]=\alpha ^n z=0$, for some $n$.
So, $\alpha=0$. Thus, $(L_1, L_0)=0$, as required.
\qed

\begin{lem}\label{[V,H]=0}
Let $L$ be a Lie superalgebra and $W$ a finite-dimensional homogeneous subspace of $L$.
If $\dim (L,L)$ is finite, then there exists a homogeneous subspace $V$ of $L$ of finite codimension such that
$(W, V)=0$.
\end{lem}
\Proof Let $w_1,\ldots, w_n$ be a homogeneous basis for $W$. Since $\dim (w_1, L)$ is finite,
there exists a homogeneous subspace $V_1$ of $L$ of finite codimension such that
$(w_1, V_1)=0$. Now we replace $L$ with $V_1$. So, there exists a homogeneous subspace $V_2$ of $V_1$ of finite codimension such that
$(w_2, V_2)=0$. Continuing this way, we get a sequence of subspaces
$$
L=V_0\supseteq V_1\supseteq V_2\supseteq\cdots \supseteq V_n,
$$
so that $\dim V_i/V_{i+1}$ is finite and   $(w_i, V_i)=0$, for every $i\geq 1$.
Clearly, $V_n$ has the desired property.
\qed

\begin{lem}\label{A,B,C}
Let $R$ be an associative algebra and $A, B, C$ some subspaces of $R$. Then
$$
[AB, _n C]\leqq \sum_{i=0}^n [A,_{i} C][B,_{n-i} C]
$$
\end{lem}
\Proof Induction on $n$.\qed

Stewart \cite{St} proved that if $H$ is a nilpotent ideal of a Lie algebra $L$ such that $L/H'$ is nilpotent then
$L$ is nilpotent. We prove the following for the Lie nilpotence of associative algebras.

\begin{prop}\label{R/I nilpotent}
Let $I$ be a two-sided ideal of an associative algebra $R$.
If $I$ is associative nilpotent of index $c$ and $R/I^2$ is Lie nilpotent of class $d$, then
$\gamma_{\mu(c, d)}(R)=0$, where $\mu(c,d)=2cd-c-d+2$.
\end{prop}
\Proof
Induction on $c$. If $c=1$, the result is obvious. So, suppose that $c\geq 2$.
Now, $I_k= I/I^{k+1}$  is an ideal of $R_k=R/I^{k+1}$.
Note that $R_k/I_k$ is Lie nilpotent and the nilpotence index of  $I_k$ is at most $c-1$, for every integer $1\leq k\leq c-1$.
 So, we may assume
$\gamma_{\mu(k, d)}(R)\su I^{k+1}$, for every $1\leq k\leq c-1$.
 By Lemma \ref{A,B,C}, we have
\begin{align*}
\gamma_{\mu(c, d)}(R)&\leqq [I^2,_{2cd -2d-c+1} R]\\
&\leqq \sum_i[I,_i R][I,_{2cd -2d-c+1-i} R],
\end{align*}
where in the sum above $i$ is in the range $0\leq i\leq 2cd -2d-c+1$.
Now we claim that every term in the sum is zero. Indeed, if  $i\leq d-1$, then
\begin{align*}
[I,_i R][I,_{2cd -2d-c+1-i} R]&\leqq I[R,_{2cd -2d-c-d+2} R]\\
&\leqq I\gamma_{\mu(c-1, d)}(R)\leqq I^{c+1}=0.
\end{align*}
On the other hand, if $i\geq d$ then
 there exists an integer $j\geq 2$ such that
$i$ belongs to the interval
$$
2(j-1)d-d-(j-1)+1\leq i< 2jd-d-j+1.
$$
Since $i\leq 2cd -2d-c+1$, we have $j\leq c$. If $j\leq c-1$ then, by induction hypothesis, we have
\begin{align*}
[I,_i R][I,_{2cd -2d-c+1-i} R]&\leqq \gamma_{\mu(j-1, d)}(R)[I,_{2(c-j)d -d-(c-j)+1} R]\\
&\leqq I^j\gamma_{\mu(c-j, d)}(R)\\
&\leqq I^jI^{c-j+1}=I^{c+1}=0.
\end{align*}
If $j= c$ then,
\begin{align*}
[I,_i R][I,_{2cd -2d-c+1-i} R]&\leqq \gamma_{\mu(c-1, d)}(R)I\\
&\leqq I^cI=I^{c+1}=0.
\end{align*}
The claim is now proved. Thus, $\gamma_{\mu(c, d)}(R)=0$, as required.
\qed

\begin{thm}\label{uL-Lie nilpotent}
Let $L=L_0\oplus L_1$ be a restricted Lie superalgebra.
Then $u(L)$ is Lie nilpotent if and only if $L'$ is finite dimensional, $L$ is a nilpotent $L_0$-module, $(L_0, L_0)$ is $p$-nilpotent,
and either $(L_1, L_1)$ is $p$-nilpotent or $\dim L_1\leq 1$  and $(L_1, L_0)=0$.
\end{thm}
\Proof
 Proposition \ref{derived-sub-f.d.}, Lemma \ref{nilp-[L,L]-p-nilp}, and Corollary  \ref{ass-PI}  are combined to give the proof of the if part.  We prove the converse. Suppose first that $\dim L_1\leq 1$  and $(L_1, L_0)=0$. Since $L_0$ is nilpotent, by induction on $\dim (L_0, L_0)$, we can assume $(L_0, L_0)=0$. It follows  that
$u(L)$ is Lie nilpotent of class at most two. Now consider the case where $(L_1, L_1)$ is $p$-nilpotent.
Let $R=u(L)$.\\

\noindent\emph{Step 1.} We may assume $L_1$ has a subspace $P$ of finite codimension such that $(L_1+(L_0, L_0), P)=0$.
 By Lemma \ref{[V,H]=0}, there exists a subspace $K$ of $L_1$ of finite codimension such that
$((L, L), K)=0$. It follows that $((K, K), L)=0$.
Let $N=\la (K, K)\ra_p\su L_0$ and set $I=NR$.
We claim it is enough to prove that $u(L/N)$ is Lie nilpotent.
Indeed, suppose $\gamma_n(u(L/N))=0$. This means that $\gamma_n(R)\su NR$.
Thus, $\gamma_{n+1}(R)\su [NR, R]=N[R, R]$, since $N$ is central in $R$.
Hence,
$$
\gamma_{nk}(R)\su \w^k(N)R.
$$
 Since $\w(N)$ is associative nilpotent, by Lemma \ref{w(L)-nilpotent}, it follows
that $R$ is Lie nilpotent.
So, we replace $L$ with $L/N$. In other words, we can assume that $(K, K)=0$.
Let $V$ be a finite-dimensional subspace of $L_1$ such that $L_1=V+K$.
By Lemma \ref{[V,H]=0}, there exists a subspace $P$ of $K$ of finite codimension such that
$(V, P)=0$. Thus, $(L_1 +(L_0, L_0), P)=0$.
Let $W$ be a finite-dimensional subspace of $L_1$ such that $L_1=W+P$.
We may assume that $(L_1, L_0)\su W$.
\\

\noindent\emph{Step 2.} We may assume $(L_1, L_1)=(L_0, L_0)=0$. 
 Consider the ideal  $N=\la W+(L_0, L_0)\ra_p$ of $L$.
By Lemma \ref{w(L)-nilpotent}, $\w(N)$ is  associative nilpotent.
Thus, $NR$, the associative ideal of $R$ generated by $N$, is nilpotent.
Furthermore, $(NR)^2=\w^2(N)R$. Now, by Proposition \ref{R/I nilpotent}, it is enough to prove that $\bar R=R/\w^2(N)R$ is Lie nilpotent. Since $\la N'\ra_p R\su \w^2(N)R$, it is enough to prove that
$R/\la N'\ra_p R=u(L/\la N'\ra_p)$ is Lie nilpotent. So, we replace $L$ with $L/\la N'\ra_p$. Thus, 
Note that  $(W, W)=(W, (L_0, L_0))=0$. It follows by Step 1 that $(L_1, L_1)=0$. 
Since $L_0$ is nilpotent, a simple induction on $\dim (L_0, L_0)$ and using Step 1 enables us to assume
that $(L_0, L_0)=0$.\\

\noindent\emph{Step 3.} $R$ is Lie nilpotent.
Note that  the subalgebra of $R$ generated by $L_1$ is the Grassmann algebra $G$ (we assume $1\in G$).
Let $A$ be the subalgebra of $R$ generated by $L_0$.
Note that $U(L)=GA$.
For every $k\geq 1$, let
$$
S_k=\sum [L_1,_{k_1} L_0]\cdots [L_1,_{k_r} L_0],
$$
 where in the  sum above the $k_i$ are all positive inters and  $\sum_{i=1}^r k_i\geq k$.
 Note that $S_kR=RS_k$, for every $k\geq 1$.
We claim that $\gamma_{2k+1}(R)\su S_{k}R$, for every $k\geq 1$.
Since $(L_1, L_0)$ is finite-dimensional and $(L_1,_n L_0)=0$, for some $n$, it follows that
 $S_k=0$, for some $k$.  So, it is enough to prove the claim. We proceed by induction on $k$.
 Since $L_0'=0$, we have
\begin{align*}
[R,R]=[GA, GA]&\su G[GA, A]+[GA, G]A\\
&\su S_1R+[G, G]A.
\end{align*}
Since $[G, G, G]=0$, we get
\begin{align*}
\gamma_3(R)&\su [S_1R, R]+[[G, G]A, GA]\\
&\su S_1R+[[G, G]A, G]A\su S_1R+[[G, G], G]A=S_1R.
\end{align*}
Now suppose that  $\gamma_{2k+1}(R)\su S_kR$. We have
\begin{align*}
[S_kGA, GA]&=G[S_kGA, A]+[S_kGA, G]A\\
&\su S_{k+1}R+[S_kG, G]A.
\end{align*}
Thus,
\begin{align*}
\gamma_{2k+3}(R)&\su  S_{k+1}R+[[S_kG, G]A, GA]\\
&\su  S_{k+1}R+G[[S_kG, G]A, A]+[[S_kG, G]A, G]A\\
&\su  S_{k+1}R+[S_kG, G][A, G]A+[[S_kG, G], G]A\\
&\su  S_{k+1}R+[[G, G], G]A=S_{k+1}R.
\end{align*}
This finishes the induction step. The proof is complete.
\qed

\end{document}